\documentclass[11pt]{article}
\usepackage{graphicx}
 \usepackage[reqno, namelimits, sumlimits]{amsmath}
   \usepackage{amsfonts}
   \usepackage{amssymb}
  \usepackage{amsthm}
   \numberwithin{equation}{section}

\newtheorem{thm}{Theorem}[section]

\newtheorem{remark}[thm]{Remark}

\def\R{{\mathbb R}}

\font\textmsbm=msbm10
\newfam\Bbbfam
\textfont\Bbbfam=\textmsbm

\date{ }
\begin{document}
\setcounter{page}{0}
\title{An oscillation result for  slightly perturbed trigonometric series.}

\maketitle

\medskip\author{ \center  Alain HARAUX \\Sorbonne Université, Université Paris Cité, \\ CNRS, Laboratoire Jacques-Louis Lions, \\LJLL, F-75005 Paris, France.\\ email: alain.haraux@sorbonne-universite.fr
 \\}

\bigskip

\renewcommand{\abstractname} {\bf Abstract}
\begin{abstract} We establish that any nontrivial anharmonic Fourier series with exponents quadratically close to the positive integers changes sign on each interval of length $2\pi$.

\end{abstract}

\bigskip

\qquad {\small \bf AMS classification numbers:}  34C10, 35B05, 35B15, 42A55, 42A75, 43A99
 
\bigskip
\qquad {\small \bf Keywords:}
oscillation theory, almost periodicity.
\vfill

\newpage

\noindent 
{\Large \bf  Introduction} \bigskip

 In the study of vibrating mechanical systems, an important role is played by almost periodic functions, cf. e.g \cite{Dieud}, section I.3 for a concise presentation. By Bochner's criterion, a real valued  function $ f $ defined on the real line 
$\R$  is almost periodic if, and only if the set of translates $f_a= f(a+.)$ is  precompact in the space of bounded continuous functions $C_b(\R, \R).$   A continuous periodic function is in particular almost periodic with a compact set of translates. The space of almost periodic functions is stable by many operations, for example a sum of any number of real valued continuous functions with arbitrary periods is almost periodic, while the sum of two non-constant periodic continuous functions with incommensurable periods is never periodic, cf. e.g.\cite{Dieud}, lemma I.3.1.7. In the case of periodic functions, it is clear that either they are constant, or they oscillate around their mean-value or average. This global property is also important in the case of general almost periodic functions.  The first case to be considered is a sum of periodic functions with mean-value 0. In \cite{HK} it was shown that a finite sum of such periodic functions has an oscillation length bounded by the sum of the periods. This result can be used to find a uniform oscillation time for an infinite sum when the sum of the periods is finite. Unfortunately, it cannot be used even in the case of vibrating strings for which a very different method had to be devised in \cite{CH1, CH2}. Another method can be applied for anharmonic Fourier series when a uniform gap separates the exponents. This was the object of \cite{osc}, relying on a technique borrowed from Ingham \cite{Ing} and improving slightly a previous one from a chapter of the lecture note \cite{H2}. We now present a slightly different approach which shows that any nontrivial anharmonic Fourier series with exponents quadratically close to the positive integers changes sign on each interval of length $2\pi$.

\section{The result }

\begin{thm} \label{main} Let $f\not \equiv 0$ be a real-valued almost periodic function with mean-value 0 defined as the limit of the uniformly convergent series
\begin{equation} \label{S}f(t) =  \sum _{n\ge 1} [ a_n \cos (\mu_nt) + b_n \sin (\mu_nt) ]. \end{equation}  We assume that 
\begin{equation}\label{eps} \sum _{n\ge 1} (\mu_n-n)^2 < \frac{1} {8\pi^2} \end{equation} Then on each interval $J$  of length  $ |J|\ge 2\pi  $,  $f$ takes both positive and negative values. \end{thm} 
\begin{proof} Since the form \eqref{S} is invariant by translation in $t$, it is sufficient to consider the case $J = [-\pi, \pi]$ .  We consider first the case  where all coefficients $b_n$ vanish, in which case $f$ is even.  This amounts to consider instead of $f$ the series \begin{equation} \label {g} g(t) =  \sum _n a_n \cos (\mu_nt) . \end{equation} We show that 
$$ g\ge 0 \,\, {\rm a.e. on}\,\, [0, \pi] \implies g \equiv 0 $$ Indeed, we deduce from the positivity condition that the function 
$$ G(t):= \int_0^t g(s) ds = \sum _n \frac{a_n}{\mu_n}  \sin (\mu_nt)$$ is non-decreasing on $[0, \pi]$, and in particular
$$ \forall t\in [0, \pi],  \quad 0 \le G(t) \le G(\pi) =  \sum _n \frac{a_n}{\mu_n} (-1) ^n  \sin (\pi \varepsilon_n) \le  \pi \sum _n |\varepsilon_n||\frac{a_n}{\mu_n}|  $$ where $$ \varepsilon_n:= \mu_n-n $$ At this stage, we observe that since the sequence $a_n$ is square summable and $\mu_n\sim n $ at infinity in $n$, the sequence $\frac{a_n}{\mu_n}$ is absolutely summable.  Let us consider the periodic function 
$$ \Gamma(t): = \sum _n \frac{a_n}{\mu_n}  \sin (n t)$$ We have, since the sine function is contractive,  
$$ \forall t\in [0, \pi],  \quad 0 \le | \Gamma(t)| \le |G(t)| +   \pi \sum _n |\varepsilon_n||\frac{a_n}{\mu_n}| \le 2\pi \sum _n |\varepsilon_n||\frac{a_n}{\mu_n}|$$ This implies 
$$ \int_0^{\pi} \Gamma(t)^2 dt  \le 4\pi^3 (\sum _n |\varepsilon_n||\frac{a_n}{\mu_n}|)^2 \le 4\pi^3 (\sum _n \varepsilon_n^2) (\sum _n |\frac{a_n}{\mu_n}|^2) $$ On the other hand , we have $$ \int_0^{\pi} \Gamma(t)^2 dt  = \frac{\pi}{2} \sum _n |\frac{a_n}{\mu_n}|^2 $$ So that we end up with the inequality 

\begin{equation} \sum _n |\frac{a_n}{\mu_n}|^2 \le 8\pi^2 \sum _n \varepsilon_n^2 \sum _n |\frac{a_n}{\mu_n}|^2\end{equation}
which together with \eqref{eps} implies  $a_n \equiv 0$. 

\bigskip\noindent  In the general case, if $f\ge 0 \,\, {\rm a.e. \,on}\,\,J = [-\pi, \pi]$, we apply the property just proved to $g(t)= f(t)+f(-t)$. It follows that $g=0$, therefore $f$ is odd and nonnegative on $J$, and we conclude that $f\equiv 0$ on $J$. 
\end{proof} 

\begin{remark} Actually, it even follows that all coefficients $a_n$ and $b_n$ vanish, since at the very last step of the proof, we end up with a series of sines which vanishes on $J$. The primitive of this series is a series of cosines with derivative 0 on J, and is therefore constant, and it particular has a constant sign. The special case of cosine series ended up with $a_n \equiv 0$ which now reads $b_n \equiv 0$.    \end{remark}

\section{Application to vibrating strings} Let us consider the simple perturbed string equation 

\begin{equation} \label{string} \begin{cases}
u_{tt}-u_{xx} -c  u = 0 \quad  \rm{on}\,\,  \R\times(0, l)   \\
u(t, 0) = u(t, l) = 0 \quad  \rm{on}\,\,  \R 
\end{cases}\end{equation} where $l>0, \,c\ge 0$. 
In this case, for any $(u_0, v_0) \in H_0^1(0, l) \times L^2(0, l) $, there is a unique solution $ u\in C(\R, H_0^1(0, l))\cap C^1(\R, L^2(0, l))\cap C^2(\R, H^{-1} (0, l))$ of \eqref{string} such that $u(0) = u_0\,\,  \rm{and}\,\, u_t(0)= v_0 $, given by the formula 
\begin{equation} u(t, x) =  \sum _{n\ge1}[ a_n \cos (\mu_nt) + b_n \sin (\mu_nt) ] \sin (\tfrac{\pi}{l} nx).\end{equation} with 
$$ u_0(x) = \sum _{n\ge1}  a_n \sin (\tfrac{\pi}{l} nx) ; \quad v_0(x) =  \sum _{n\ge1} \mu_n  b_n  \sin (\tfrac{\pi}{l} nx) $$ and $$ \mu_n := \sqrt{(\tfrac{\pi}{l})^2n^2-c}$$ The application of Theorem \ref{main} implies  that if $c$ is small enough, whenever $\frac{x}{l}$ is irrational and $u\not\equiv 0$, the function $t\longrightarrow u(t, x)$ takes both positive and negative values on the real interval $J$ as soon as $ |J|\ge 2l$. 

\begin{remark}   It is not sufficient to assume $c< (\tfrac{\pi}{l})^2. $  For instance for $l = \pi$, the function $ u(t, x) =  \sin \frac{t}{2} \sin x $  is a solution of \eqref {string} with $c = \frac{3}{4} $ with $ u(t, x) \ge 0 $ on $[0, 2\pi]  \times [0, \pi]. $  \end{remark}

\begin{remark}    A sufficient condition for $l = \pi$ is $$ c^2 \sum _1^{\infty} \frac {1} { (n+ \sqrt{n^2-c})^2} < \frac{1}{8\pi^2}  $$ In particular the condition 
$ c< \frac{\sqrt{3}}{2\pi^2}\sim 0.0877 $ is obviously sufficient but quite strong. Numerically, it can be checked that $ c<0.17$ is sufficient.  It would be interesting to find the best upper bound for c. \end{remark}


\begin{thebibliography}{10}


\bibitem{AP} L. Amerio and G. Prouse, \emph{Almost-periodic functions and functional equations}, Van Nostrand Reinhold Co., New York-Toronto, Ont.-Melbourne 1971 i-viii+184 pp. 

\bibitem{Bes} A.S. Besicovitch, \emph{ Almost Periodic Functions}; Dover: New York, NY, USA, 1954. 

\bibitem{BVN} S. Bochner and J. Von Neumann, \emph{On compact solutions of operational-differential equations}, Ann. of Math. \textbf{2}, 36 (1935), 255--291.

\bibitem{CH1} T. Cazenave and A. Haraux, \emph{Propri{\' e}t{\' e}s oscillatoires des solutions de certaines {\' e}quations des ondes semi-lin{\' e}aires}, C. R. Acad. Sci. Paris Ser. I Math. \textbf{298} (1984), no. 18, 449--452.

\bibitem{CH2} T. Cazenave and A. Haraux, \emph{Oscillatory phenomena associated to semilinear wave equations in one spatial dimension}, Trans. Amer. Math. Soc. \textbf{300} (1987), no. 1, 207--233.

\bibitem{Cord1} C.Corduneanu, \emph{Almost periodic functions}, Interscience Tracts in Pure and Applied Mathematics \textbf{22},  [John Wiley \& Sons], New-York-London-Sydney (1968).

\bibitem{Fink} A.M. Fink, \emph{Almost periodic differential equations}, Lecture Notes in Math., Vol. 377, Springer-Verlag, Berlin-New York, 1974. i-viii+336 pp.

\bibitem {H1} A. Haraux,  \emph{A simple almost-periodicity criterion and applications}, {\it J. Diff. Eq. } {\bf 66} (1987), no.1, 51--61. 

\bibitem {Dieud} A. Haraux,  \emph{Semi-linear hyperbolic problems in bounded domains},  Math. Rep. 3 (1987), no. 1, i--xxiv and 1--281. 

\bibitem {H2} A. Haraux,  \emph {Quelques propriétés des séries lacunaires utiles dans l'étude des vibrations élastiques} in "Nonlinear partial differential equations and their applications, Coll\`{e}ge de France Seminar volume XII, Pitman research Notes in Mathematics N°302 (1994), 113--124.

\bibitem {OscPort} A. Haraux, \emph{On the strong oscillatory behavior of all solutions to some second order evolution equations}, Port. Math. {\bf 72} (2015), no. 2, 193--206. 

\bibitem {RioLNt} A. Haraux, \emph{Nonlinear vibrations and the wave equation}, SpringerBriefs Math.BCAM SpringerBriefs, Springer, Cham, (2018) i-x+102 pp.

\bibitem {Op} A. Haraux,  \emph {Some simple problems for the next generations}, Partial differential equations arising from physics and geometry, 296–310.
London Math. Soc. Lecture Note Ser. {\textbf 450}, Cambridge University Press, Cambridge, 2019.

\bibitem {osc} A. Haraux,  \emph {Uniform oscillation properties of almost periodic functions}, arXiv:2610.06515 [math.AP], 2026.

 \bibitem{HK} A. Haraux  and V. Komornik, \emph{ Oscillations of anharmonic Fourier series and the wave equation},  Rev. Mat. Iberoamericana \textbf{1} (1985), no. 4, 57--77.
 
 \bibitem{Ing} A.E. Ingham, \emph{ Some trigonometrical inequalities with applications to the theory of series}, Math Z 41 (1936), 367–379 .
 
 \bibitem{Yosh} T. Yoshizawa, \emph{Stability theory and the existence of periodic solutions and almost periodic solutions},  Appl. Math. Sci., Vol. 14
Springer-Verlag, New York-Heidelberg, 1975. i-vii+233 pp.

\end{thebibliography}
\end{document}